\documentclass{amsart}
\usepackage{latexsym,amsmath,amssymb}

\theoremstyle{definition}
\theoremstyle{remark}

\numberwithin{equation}{section}

\begin{document}

\title[COMPACT LAMBERT TYPE OPERATORS]
{COMPACT LAMBERT TYPE OPERATORS BETWEEN TWO $L^{p}$ SPACES}

\author{\sc Y. Estaremi and M. R. Jabbarzadeh }
\address{\sc y. estaremi and  m. r. jabbarzadeh}
\email{estaremi@gmail.com} \email{mjabbar@tabrizu.ac.ir}

\address{ faculty of mathematical
sciences, university of tabriz, p. o. box: 5166615648, tabriz, iran}

\thanks{}

\thanks{}



\date{}

\dedicatory{}

\commby{}

\begin{abstract}
In this paper we investigate compactness of weighted
Lambert type operators  between two $L^{p}$ spaces.

\noindent {}
\end{abstract}

\maketitle

\section{ \sc\bf Introduction and Preliminaries}

Let $(X,\Sigma,\mu)$ be a complete $\sigma$-finite measure space.
For any sub-$\sigma$-finite algebra $\mathcal{A}\subseteq
 \Sigma$ with $1\leq p\leq \infty$, the $L^p$-space
$L^p(X,\mathcal{A},\mu|_{\mathcal{A}})$ is abbreviated  by
$L^p(\mathcal{A})$, and its norm is denoted by $\|.\|_p$. All
comparisons between two functions or two sets are to be
interpreted as holding up to a $\mu$-null set. The support of a
measurable function $f$ is defined as $\sigma(f)=\{x\in X;
f(x)\neq 0\}$. We denote the
vector space of all equivalence classes of almost everywhere finite
valued measurable
functions on $X$ by $L^0(\Sigma)$.

\vspace*{0.3cm} For a sub-$\sigma$-finite algebra
$\mathcal{A}\subseteq\Sigma$, the conditional expectation
operator associated with $\mathcal{A}$ is the mapping
$f\rightarrow E^{\mathcal{A}}f$, defined for all non-negative $f$
as well as for all $f\in L^p(\Sigma)$, $1\leq p\leq \infty$,
where $E^{\mathcal{A}}f$, by the Radon-Nikodym theorem, is the
unique $\mathcal{A}$-measurable function satisfying
$$\int_{A}fd\mu=\int_{A}E^{\mathcal{A}}fd\mu, \ \ \ \forall A\in \mathcal{A} .$$
As an operator on $L^{p}({\Sigma})$, $E^{\mathcal{A}}$ is an
idempotent and $E^{\mathcal{A}}(L^p(\Sigma))=L^p(\mathcal{A})$.
If there is no possibility of confusion we write $E(f)$ in place
of $E^{\mathcal{A}}(f)$. This operator will play major role in
our work, and we list here some of its useful properties:

\vspace*{0.2cm} \noindent $\bullet$ \  If $g$ is
$\mathcal{A}$-measurable then $E(fg)=E(f)g$.

\noindent $\bullet$ \ $|E(f)|^p\leq E(|f|^p)$.

\noindent $\bullet$ \ If $f\geq 0$ then $E(f)\geq 0$; if $f>0$
then $E(f)>0$.

\noindent $\bullet$ \ $|E(fg)|\leq
E(|f|^p)|^{\frac{1}{p}}E(|g|^{p'})|^{\frac{1}{p'}}$, where
$\frac{1}{p}+\frac{1}{p'}=1$ H\"{o}lder inequality).

\noindent $\bullet$ \ For each $f\geq 0$, $\sigma(f)\subseteq
\sigma(E(f))$.

\vspace*{0.2cm}\noindent A detailed discussion and verification
of most of these properties may be found in \cite{rao}. We recall
that an $\mathcal{A}$-atom of the measure $\mu$ is an element
$A\in\mathcal{A}$ with $\mu(A)>0$ such that for each
$F\in\mathcal{A}$, if $F\subseteq A$ then either $\mu(F)=0$ or
$\mu(F)=\mu(A)$. A measure space $(X,\Sigma,\mu)$ with no atoms is
called non-atomic measure space. It is well-known fact that every
$\sigma$-finite measure space $(X,
\mathcal{A},\mu_{\mid_{\mathcal{A}}})$ can be partitioned
uniquely as $X=\left (\bigcup_{n\in\mathbb{N}}A_n\right )\cup B$,
where $\{A_n\}_{n\in\mathbb{N}}$ is a countable collection of
pairwise disjoint $\mathcal{A}$-atoms and $B$, being disjoint
from each $A_n$, is non-atomic (see \cite{z}).

\vspace*{0.2cm} Let $u\in L^0(\Sigma)$. Then $u$ is said to be
conditionable with respect to $E$ if $u\in\mathcal{D}(E):=\{f\in
L^0(\Sigma): E(|f|)\in L^0(\mathcal{A})\}$.
 Take $u$ and $w$ in $\mathcal{D}(E)$. As it is defined in \cite{es}, weighted lambert type operator $T$
 is bounded operator $T:=M_{w}EM_{u}$ from $L^p(\Sigma)$ into $L^q(\Sigma)$,
 where $M_u$ and $M_w$ are multiplication operators.
Throughout this paper we assume that $u$ and $w$ are in $\mathcal{D}(E)$, $E=E^\mathcal{A}$
and $T$ is weighted lambert type operator.

\vspace*{0.2cm} Combination of conditional expectation operator
$E$ and multiplication operators  appears more often in the
service of the study of other operators such as multiplication
operators and weighted composition operators. Specifically, in
\cite{mo}, S.-T. C. Moy has characterized all operators on $L^p$
of the form $f\rightarrow E(fg)$ for $g$ in $L^q$ with $E(|g|)$
bounded. In \cite{dou}, R. G. Douglas analyzed positive
projections on $L\sb {1}$ and many of his characterizations are
in terms of combinations of multiplications and conditional
expectations.

\vspace*{0.2cm} Some classic properties of the operator $EM_{u}$ on $L^p(\Sigma)$ spaces is characterized in \cite{he, her} and \cite{lam}. The authors have characterized boundedness of $T=M_{w}EM_{u}$ between two $L^p(\Sigma)$ spaces, polar decomposition of $T=M_{w}EM_{u}$ on $L^2(\Sigma)$ and some other classic properties of $T=M_{w}EM_{u}$ on $L^p(\Sigma)$ \cite{es}.

In this paper we characterize compactness of $T=M_{w}EM_{u}$ between two $L^p(\Sigma)$ spaces.

\section{ \sc\bf  Compact Weighted
Lambert Type Operators }

\vspace*{0.3cm} {\bf Theorem 2.1.} Let $1<q<p<\infty$ and let $p',
q'$ be conjugate component to $p$ and $q$ respectively. Then
weighted Lambert type operator $T=M_wEM_u$ from $L^{p}(\Sigma)$ into
$L^{q}(\Sigma)$ is compact if and only if\\

(i) $(E|w|^{q})^{{1}/{q}}(E|u|^{p'})^{{1}/{p'}}=0$ a.e. on $B$;\\

(ii) $\sum_{n\in\mathbb{N}}(E(|w|^{q}))^{{p'q'}/{(q'-p')}}(A_{n})
(E(|u|^{p'}))^{{q'}/{(q'-p')}}(A_{n})\mu(A_{n})<\infty$.

\vspace*{0.3cm} {\bf Proof.} Let $f\in L^{p}(\Sigma)$. Then
$$\|Tf\|^{q}_{q}=\int_{X}|w|^{q}|E(uf)|^{q}d\mu=\int_{X}E(|w|^{q})|E(uf)|^{q}d\mu$$$$=
\int_{X}|E(u(E(|w|^{q}))^{\frac{1}{q}} f)|^{q}d\mu
=\|EM_{v}f\|^{q}_{q},$$ where $v:=u(E(|w|^{q}))^{{1}/{q}}$. It
follows that the weighted conditional operator $T: L^{p}(\Sigma)\rightarrow L^{q}(\Sigma)$
is compact if and only if
$R_v^*=M_{\bar{v}}E=M_{\bar{v}}:L^{q'}(\mathcal{A})\rightarrow
L^{p'}(\Sigma)$ is compact. Also, since
$\|M_{\bar\nu}f\|_{p'}=\|M_{(E(|v|^{p'}))^{1/p'}}f\|_{p'}$,
$M_{\bar{v}}$ is compact if and only if the multiplication operator $M_{(E(|v|^{p'}))^{1/p'}}$ is compact.
Now, suppose that $T$ is bounded.
Firstly, we show that
$E(|v|^{p'})=0$ a.e on $B$. Suppose, on the contrary. Then, there
exists some $\delta>0$ such that the set $B_{0}=\{x\in
B:E(|v|^{p'})(x)>\delta\}$ has positive measure. We
may also assume $\mu(B_{0})<\infty$.
Since $B_0$ has no $\mathcal A$-atoms, we can find $E_n\in{\mathcal A}$
of positive measure satisfying $E_{n+1}\subseteq E_n\subseteq B_0$
for all $n$ and $\lim_n\mu(E_n)=0$.
Put $f_n=\chi_{E_n}/(\mu(E_n))^{1/p'}$. Then for each $n$, $f_n$ is bounded element in $L^{q'}(\mathcal A)$ and $\lim_nf_n(x)=0$ for all $x\in X\setminus\cap_nE_n$.
It follows that $lim_nE(|v|^{p'})f_n=0$ pointwise almost everywhere, because $\mu(\cap_nE_n)=0$.
By compactness of $M_{(E(|v|^{p'}))^{1/p'}}$, $\{(E(|v|^{p'}))^{1/p'}f_n\}_n$ is uniformly integrable in $L^{p'}(\Sigma)$. Hence there exists $t\geq 0$ such that
$$\sup_{n\in\mathbb{N}}\int_{\{E(|v|^{p'})f_n\geq\ t\}}E(|v|^{p'})f_n^{p'}d\mu<\delta,$$
and so $\int_{E_N}E(|v|^{p'})f_n^{p'}d\mu<\delta$ for some large $N\in\mathbb{N}$. Then we get that
$$\delta>\int_{E_N}E(|v|^{p'})f_n^{p'}d\mu=\|M_{(E(|v|^{p'}))^{1/p'}}f_n^{p'}\|^{p'}=\frac{1}{\mu(E_N)}
\int_{E_N}E(|v|^{p'})\geq\delta ,$$
but this is a contradiction. Next, we examine (ii). Since $T$ is compact
then $M_{v}$ is bounded, and so by Theorem 2.2 in \cite{es},
$E(|v|^{p'})^{{1}/{p'}}\in L^{{p'q'}/{(q'-p')}}(\mathcal{A})$. Thus $\sum_{n\in\mathbb{N}}(E(|v|^{p'}))^{{q'}/{(q'-p')}}(A_{n})\mu(A_{n})<\infty$.

\vspace*{0.3cm}
Conversely, assume both (i) and (ii) hold.
By the same preceding discission, it suffices to establish the
compactness of
$M_{(E(|v|^{p'}))^{{1}/{p'}}}:L^{q'}(\mathcal{A})\rightarrow
L^{p'}(\Sigma)$. From (ii), for any $\varepsilon>0$, there exists
$N_{\varepsilon}\in\mathbb{N}$ such that
$\sum_{n>N_{\varepsilon}}(E(|v|^{p'}))^{{q'}/{(q'-p')}}(A_{n})\mu(A_{n})<\varepsilon.$
Put $v_{\varepsilon}:=\sum_{n\leq
N_{\varepsilon}}(E(|v|^{p'}))^{{1}/{p'}}(A_{n})\chi_{A_{n}}$.
Obviously, $M_{v_{\varepsilon}}$ is a bounded and finite rank operator
from $L^{q'}(\mathcal{A})$ to $L^{p'}(\Sigma)$. Moreover, by H\"{o}lder
inequality, for each $f\in L^{q'}(\mathcal{A})$ we have
$$\|M_{(E(|v|^{p'}))^{{1}/{p'}}}f-M_{v_{\varepsilon}}f\|^{p'}_{p'}=
\int_{\cup_{n>N_{\varepsilon}}A_{n}}E(|v|^{p'})|f|^{p'}d\mu$$$$=
\sum_{n>N_{\varepsilon}}E(|v|^{p'})(A_{n})\mu(A_{n})|f(A_{n})|^{p'}=
\sum_{n>N_{\varepsilon}}E(|v|^{p'})(A_{n})\mu(A_{n})^{\frac{q'-p'}{q'}}|f(A_{n})|^{p'}\mu(A_{n})^{\frac{p'}{q'}}$$$$
\leq(\sum_{n>N_{\varepsilon}}E(|v|^{p'})^{\frac{q'}{q'-p'}}(A_{n})\mu(A_{n}))^{\frac{q'-p'}{q'}}
(\sum_{n>N_{\varepsilon}}|f(A_{n})|^{q'}\mu(A_{n}))^{\frac{p'}{q'}}\leq\varepsilon^{\frac{q'-p'}{q'}}\|f\|^{p'}_{q'}.$$
It follows that
$\|M_{(E(|v|^{p'}))^{{1}/{p'}}}-M_{v_{\varepsilon}}\|\leq\varepsilon^{{(q'-p')}/{q'}}$.
So $M_{(E(|v|^{p'}))^{{1}/{p'}}}$ is the limit of some
finite rank operators and is therefore compact.
This completes the proof of the theorem.\ $\Box$

\vspace*{0.3cm} {\bf Theorem 2.2.} Let $1< p< q<\infty$ and let $p',
q'$ be conjugate component to $p$ and $q$ respectively. Then the
 weighted Lambert type operator $T$ from
$L^{p}(\Sigma)$ into $L^{q}(\Sigma)$ is compact if and only if\\

(i) $E(|u|^{p'})^{{1}/{p'}}(E(|w|^{q}))^{{1}/{q}}=0$ a.e. on
$B$;\\

(ii) $\lim_{n\rightarrow\infty}\frac{E(|u|^{p'})(A_{n})(E(|w|^{q}))^{{p'}/{q}}(A_{n})}
{\mu(A_{n})^{{(p'-q')}/{q'}}}=0 $, when the number of $\mathcal{A}$-atoms is not finite.

\vspace*{0.3cm} {\bf Proof.} Let $f\in L^{p}(\Sigma)$.
By the same argument in the proof of Theorem 2.1, since $\|Tf\|_{q}=\|EM_{v}f\|_{q}$
$v:=u(E(|w|^{q}))^{\frac{1}{q}}$, it follows that
$T:L^{p}(\Sigma)\rightarrow L^{q}(\Sigma)$
is compact if and only if
$M_{v}:L^{q'}(\mathcal{A})\rightarrow
L^{p'}(\Sigma)$ is compact. Let $T$ is compact. Hence $M_{v}$ is compact and so is bounded.
Thus by Theorem 2.3 in \cite{es}, $E(|v|^{p'})=0$ a.e on $B$ and thus (i) is holds.
Now, if (ii) is not hold,
we can fined a constant $\delta>0$ such that
${E(|v|^{p'})(A_{n})}/{\mu(A_{n})^{{(p'-q')}/{q'}}}>\delta$,
 for all $n\in\mathbb{N}$. For each $n\in\mathbb{N}$, define
 $f_{n}={\chi_{A_{n}}}/{\mu(A_{n})^{{1}/{q'}}}$. Obviously,
 $f_{n}\in L^{q'}(\mathcal{A})$ and $\|f_{n}\|_{q'}=1$. Thus, for each $m,n\in\mathbb{N}$ with $m\neq n$ we
get that
$$\|M_{\bar{v}}f_{m}-M_{\bar{v}}f_{n}\|^{p'}_{p'}=\int_{X}|v|^{p'}|f_{m}-f_{n}|^{p'}d\mu=
\int_{X}E(|v|^{p'})|f_{m}-f_{n}|^{p'}d\mu$$
$$\geq\int_{A_{m}}E(|v|^{p'})|f_{m}|^{p'}d\mu+\int_{A_{n}}E(|v|^{p'})|f_{n}|^{p'}d\mu\geq2\delta.$$
But this is a contradiction.

\vspace*{0.3cm}
Conversely, assume both (i) and (ii) hold. Since
By the same argument in the proof of Theorem 2.1, it suffices to establish the
compactness of
$M_{(E(|v|^{p'}))^{{1}/{p'}}}$. From (ii), for any $\varepsilon>0$, there exists
$N_{\varepsilon}\in\mathbb{N}$ such that
${E(|v|^{p'})(A_{n})}/{\mu(A_{n})^{{(p'-q')}/{q'}}}<\varepsilon$,
for all $n>N_{\varepsilon}$.
Take $v_{\varepsilon}:=\sum_{i\leq
N_{\varepsilon}}E(|v|^{p'})^{{1}/{p'}}(A_{i})\chi_{A_{i}}$.
Then, $M_{v_{\varepsilon}}$ is a bounded and finite rank operator
from $L^{q'}(\mathcal{A})$ to $L^{p'}(\Sigma)$. Moreover, for each
$f\in L^{q'}(\mathcal{A})$ with $\|f\|_{q'}=1$, we have
$|f(A_{n})|^{q'}\mu(A_{n})\leq\|f\|^{q'}_{q'}=1$.
 Since ${p'}/{q'}>1$, then
$(|f(A_{n})|^{q'}\mu(A_{n}))^{{p'}/{q'}}\leq|f(A_{n})|^{q'}\mu(A_{n})$.
Thus we conclude that
$$\|M_{(E(|v|^{p'}))^{{1}/{p'}}}f-M_{v_{\varepsilon}}f\|^{p'}_{p'}=
\int_{\cup_{n>N_{\varepsilon}}A_{n}}E(|v|^{p'})|f|^{p'}d\mu$$$$=\sum_{n>N_{\varepsilon}}
E(|v|^{p'})(A_{n})\mu(A_{n})|f(A_{n})|^{p'}
$$$$=\sum_{n>N_{\varepsilon}}\frac{E(|v|^{p'})(A_{n})}{\mu(A_{n})^{{(p'-q')}/{q'}}}
(|f(A_{n})|^{q'}\mu(A_{n}))^{\frac{p'}{q'}}
\leq\varepsilon\|f\|^{p'}_{q'}.$$
Then
$\|M_{(E(|v|^{p'}))^{{1}/{p'}}}-M_{v_{\varepsilon}}\|\leq\varepsilon^{{1}/{q'}}$,
and so $M_{(E(|v|^{p'}))^{{1}/{p'}}}$ is compact.\ $\Box$

\vspace*{0.3cm}
{\bf Remark 2.3.}
(a) According to the procedure used in the proof of Theorem 2.1,
the weighted Lambert type operator $T=M_wEM_u$ from $L^{p}(\Sigma)$ into
$L^{p}(\Sigma)$ is compact if and only if  $\nu=0$ a.e. on $B$ and  $\lim_{n\rightarrow\infty}\nu(A_{n})=0$,
where $\nu=(E|w|^{p})^{{1}/{p}}(E|u|^{p'})^{{1}/{p'}}$. See Lemma 2.5 in \cite{es} for another characterization.\\

(b) If $q=1$, then
by according Theorem 2.2, the
weighted Lambert type operator $T$ from
$L^{p}(\Sigma)$ into $L^{1}(\Sigma)$ is compact if and only if
$E(|u|^{p'})^{{1}/{p'}}E(|w|)=0$ a.e. on $B$ and
$\sum_{n\in\mathbb{N}}E(|u|^{p'})(A_{n})(E(|w|))^{p'}(A_{n})\mu(A_{n})<\infty$.\\

(c) Let $(X, \mathcal{A}, \mu)$ be a non-atomic measure space. Since for each $\alpha>0$ and $\beta>0$,
$\sigma(w)=\sigma(|w|^\alpha)\subseteq\sigma(E(|w|^\alpha))=\sigma((E(|w|^\alpha))^\beta)$, Then by the pervious
results, the weighted Lambert type operator $T=M_wEM_u$ from $L^{q}(\Sigma)$ into
$L^{p}(\Sigma)$ with $1<p<\infty$ and $1\leq q<\infty$ is compact if and only if it is a zero operator.

\vspace*{0.3cm} {\bf Theorem 2.4.}  Let $1<q<\infty$ and
let $X=\left (\cup_{n\in\mathbb{N}}C_n\right )\cup C$, where
$\{C_n\}_{n\in\mathbb{N}}$ is a countable collection of pairwise
disjoint $\Sigma$-atoms and $C\in \Sigma$, being disjoint from
each $C_n$, is non-atomic.
If the
 weighted Lambert type operator $T$ from
$L^{1}(\Sigma)$ into $L^{q}(\Sigma)$ is compact then\\

(i) $E(|u|^{q'})^{{1}/{q'}}(E(|w|^{q}))^{{1}/{q}}=0$ a.e. on
$B$;\\

(ii) $\lim_{n\rightarrow\infty}\frac{E(|u|^{q'})(A_{n})(E(|w|^{q}))
^{{q'}/{q}}(A_{n})}{\mu(A_{n})}=0 $, when the number of $\mathcal{A}$-atoms is not finite.
On the other hand, if\\

(i) $u(E(|w|^{q}))^{{1}/{q}}=0$ a.e. on
$C$, and\\

(ii) $\lim_{n\rightarrow\infty}\frac{E(|u|^{q'})(C_{n})(E(|w|^{q}))^{{q'}/{q}}
(C_{n})}{\mu(C_{n})}=0 $, when the number of $\Sigma$-atoms is not finite,
then weighted Lambert type operator $T$ from $L^{1}(\Sigma)$ into
$L^{q}(\Sigma)$ is compact.

\vspace*{0.3cm} {\bf Proof.} (a) Let $f\in L^{1}(\Sigma)$. Then
$\|Tf\|_{q}=\|EM_{v}f\|_{q}$,
where $v:=u(E(|w|^q))^{\frac{1}{q}}$. It follows that $T:L^{1}(\Sigma)\rightarrow L^{q}(\Sigma)$
is compact if and only if the operator $EM_{v}:L^{1}(\Sigma)\rightarrow L^{q}(\Sigma)$ is compact. Suppose that $EM_{v}$ is compact.
Thus by Theorem 2.4(b) in \cite{es}, $E(|v|^{q'})=0$ a.e on $B$. Now, we examine (ii).
Assume on the contrary, thus  we can fined a
constant $\delta>0$ such that
 ${E(|v|^{p'})(A_{n})}/{\mu(A_{n})}>\delta$,
 for all $n\in\mathbb{N}$. For every $n\in\mathbb{N}$, define
 $f_{n}={\chi_{A_{n}}}/{\mu(A_{n})^{{1}/{q'}}}$. Then
 $f_{n}\in L^{q'}(\mathcal{A})$ and $\|f_{n}\|_{q'}=1$. Thus, for any $m,n\in\mathbb{N}$ with $m\neq n$,
we obtain
$$\|M_{\bar{v}}f_{m}-M_{\bar{v}}f_{n}\|^{q'}_{\infty}
\geq\frac{1}{\mu(A_{n})^2}\int_{A_{n}}E(|v|^{p'})d\mu=\frac{E(|v|^{p'})(A_{n})}{\mu(A_{n})}\geq\delta.$$
But this is a contradiction.

\vspace*{0.3cm}
(b) Assume both (i) and (ii) hold. Then for any $\varepsilon>0$,
there exists some $N_{\varepsilon}\in \mathbb{N}$ such that
${E(|v|^{q'})(C_{n})}/{\mu(C_{n})}<\varepsilon.$
Take $v_{\varepsilon}:=\sum_{n\leq N_{\varepsilon}}v(C_{n})\chi_{C_{n}}$.
It follows that $EM_{v_{\varepsilon}}$ is a bounded and finite rank operator.
Moreover, for any $f\in L^{1}(\Sigma)$ with $\|f\|_{1}=1$  we get that
$|f(A_{n})|\mu(A_{n})\leq\|f\|_{1}=1$.
Since $q>1$, then
$(|f(A_{n})|\mu(A_{n}))^{q}\leq|f(A_{n})|\mu(A_{n})$.
Now, by conditional type H\"{o}lder inequality, we have
$$\|EM_{v}f-EM_{v_{\varepsilon}}f\|^{q}_{q}=\int_{X}|E(v\chi_{\cup_{n>N_{\varepsilon}}C_{n}}f)|^qd\mu$$
$$\leq\int_{X}(E(|v|^{q'}))^{\frac{q}{q'}}E(\chi_{\cup_{n>N_{\varepsilon}}C_{n}}|f|^q)d\mu=
\int_{\cup_{n>N_{\varepsilon}}C_{n}}(E(|v|^{q'}))^{\frac{q}{q'}}|f|^qd\mu$$
$$\sum_{n>N_{\varepsilon}}\frac{(E(|v|^{q'}))^{\frac{q}{q'}}(C_{n})}{\mu(C_{n})^{q-1}}(|f(C_{n})|\mu(C_{n}))^q\leq
\varepsilon^{\frac{q}{q'}}\sum_{n>N_{\varepsilon}}|f(C_{n})|\mu(C_{n})\leq\varepsilon^{\frac{q}{q'}}.$$
This implies that $EM_{v}$ is compact.\ $\Box$

\vspace*{0.3cm} {\bf Example  2.5.}
(a) Let $X=[0,\infty)\times
[0,\infty)$, $d\mu=dxdy$, $\Sigma$  the  Lebesgue subsets of $X$ and let
$\mathcal{A}=\{A\times [0,\infty): A \ \mbox{is a Lebesgue set in} \
[0,\infty)\}$. Put $w=1/y$ and $u=1$. Then, for each $f$ in $L^p(\Sigma)$ with $p>1$,
$$(Tf)(x, y)=\frac{1}{y}E(f)(x, y)=\frac{1}{y}\int_0^\infty f(x, t)dt .$$
Since $\int_0^\infty 1/t^pdt=E(|w|^p)\notin L^\infty(\mathcal{A})$,
so the averaging operator $T$ is not compact because
it is not bounded (see \cite[Theorem 2.1(a)]{es}).\\

(b) Let $\nu=\{m_n\}_{n=1}^{\infty}$
be a sequence of positive real numbers. Consider the measure space
$(\mathbb{N}, 2^{\mathbb{N}},\mu)$ with $\mu(\{n\})=m_n$.
Let $\varphi:\mathbb{N}\rightarrow \mathbb{N}$ be a
measurable transformation. Take $\mathcal{A}=\varphi^{-1}(2^{\mathbb{N}})$. Then for all
non-negative sequence $f=\{f_n\}_{n=1}^\infty$ and $k\in\mathbb{N}$, we have
$$  E(f)(k)=\frac{\sum_{j\in \varphi^{-1}(\varphi
(k))}f_jm_j}{\sum_{j\in \varphi^{-1}(\varphi (k))}m_j}\ . $$
Now, set $\varphi=\chi_{\{1\}}+(n-1)\chi_{\{n\in\mathbb{N}: n\geq2\}}$ and let $\mu(\{n\})=1$
for all $n\in\mathbb{N}$. Then $\mathcal{A}$ is generated by the atoms
$$\{1, 2\}, \{3\}, \{4\}, \ \ldots .$$
Hence for $w=\{w_n\}\in l^0$, $u=\{u_n\}\in\mathcal{D}(E)$ and $f=\{f_n\}\in L^p(\Sigma)$ we have
$$Tf=wE(uf)=(\frac{w_1(u_1f_1+u_2f_2)}{2}, \frac{w_2(u_1f_1+u_2f_2)}{2}, w_3u_3f_3, w_4u_4f_4, \ldots ) .$$
In particular, by Theorem 2.1, $T: L^3(\Sigma)\rightarrow L^2(\Sigma)$ is compact if and only if
$\sum_{n=3}^\infty(w_nu_n)^6<\infty$. Whence, by Theorem 2.2, $T: L^2(\Sigma\rightarrow L^3(\Sigma)$
is compact if and only if $w_nu_n\rightarrow 0$ as $n\rightarrow\infty$. \\

(c) Let $X=\mathbb{N}$,
$\Sigma=2^{\mathbb{N}}$ and let $\mu$ be the counting measure. Put
$$A=\{\{2\},\{4,6\},
\{8,10,12\},\{14,16,18,20\},\cdots\}\cup\{\{1\}, \{3\}, \{5\},
\cdots \}.$$
If we let $A_{1}=\{2\}$, $A_{2}=\{4,6\}$,
$A_{3}=\{8,10,12\}$, $\cdots$, then we see that $\mu(A_{n})=n$
and for every $n\in \mathbb{N}$, there exists $k_{n}\in
\mathbb{N}$ such that $A_{n}=\{2k_{n}, 2(k_{n}+1), \cdots,
2(k_{n}+n-1)\}$. Let $\mathcal{A}$ be the $\sigma$-algebra
generated by the partition $A$ of $\mathbb{N}$. Note that,
$\mathcal{A}$ is a sub-$\sigma$-finite algebra of $\Sigma$ and
each of element of $\mathcal{A}$ is an $\mathcal{A}$-atom. It is
known that the conditional expectation of any  $f\in
\mathcal{D}(E)$ relative to $\mathcal{A}$ is
$$E(f)=\sum_{n=1}^{\infty}\left(\frac{1}{\mu(A_n)}\int_{A_n}fd\mu\right)\chi_{A_n}+
\sum_{n=1}^{\infty}f(2n-1)\chi_{\{2n-1\}}.$$
Define $u(n)=n$ and $w(n)=\frac{1}{n^{3}}$, for all $n\in \mathbb{N}$.
For each even number $m\in \mathbb{N}$, there exists $n_{m}\in \mathbb{N}$ such that
$m\in A_{n_{m}}$. Thus for all $1\leq p,q<\infty$ we get that
$$E(|w|^q)(m)=\frac{1}{(2k_{n_m})^{3q}}+\cdots
+\frac{1}{(2k_{n_m}+2n_m-2)^{3q}}$$
and
$$E(|u|^{p})(m)=2^pk^p_{n_m}+2^p(k_{n_m}+1)^p+ \cdots+
2^p(k_{n_m}+n_m-1)^p$$
Since $n_m\leq k_{n_m}$ we have

$$(E(|w|^q))^{\frac{1}{q}}(m)(E(|u|^{p}))^{\frac{1}{p}}(m)\leq\frac{4k_{n_m}}{2^3k^3_{n_m}}.$$

Also, for all $n\in \mathbb{N}$
$$(E(|w|^q))^{\frac{1}{q}}(2n-1)(E(|u|^{p}))^{\frac{1}{p}}(2n-1)\leq
\frac{1}{(2n-1)^2} .$$
It is easy to see that
$\lim_{n\rightarrow\infty}\frac{E(|u|^{p'})(A_{n})(E(|w|^{q}))^{{p'}/{q}}(A_{n})}
{n^{{(p'-q')}/{q'}}}=0 $
and
$\lim_{n\rightarrow\infty}E(|u|^{p'})(2n-1)(E(|w|^{q}))^{{p'}/{q}}(2n-1)=0$.

Thus by theorem 2.2 the operator
$T=M_wEM_u$ is a compact weighted Lambert type operator from
$L^p(\Sigma)$ into $L^q(\Sigma)$, when $1<p<q<\infty$.
Also, we have
$\lim_{n\rightarrow\infty}E(|u|^{q'})(n)(E(|w|^{q}))^{{q'}/{q}}(n)=0$. By theorem 2.4 the operator $T=M_wEM_u$ is a compact operator from $L^{1}(\Sigma)$ into
$L^{q}(\Sigma)$, for $1<q<\infty$.


\vspace*{0.3cm}

\begin{thebibliography}{99}

\bibitem{dou}
 R. G. Douglas, Contractive projections on an $L\sb{1}$ space,
 Pacific J. Math. {\bf 15} (1965), 443-462.

\bibitem{es}
Y. Estaremi and M. R. Jabbarzadeh,
Weighted Lambert type operators on $L^{p}$ spaces,
to appear in Oper. Matrices.

\bibitem{he}
J. Herron, Weighted conditional expectation operators on $L^p$
spaces, UNC Charlotte Doctoral Dissertation, 2004.

\bibitem{her}
J. Herron, Weighted conditional expectation operators, Oper. Matrices {\bf 1} (2011), 107-118.

\bibitem{lam}
A. Lambert, $L^p$ multipliers and nested sigma-algebras, Oper.
Theory Adv. Appl. {\bf 104} (1998),  147-153.


\bibitem{mo}
Shu-Teh Chen, Moy,  Characterizations of conditional expectation
as a transformation on function spaces,
 Pacific J. Math. {\bf 4} (1954), 47-63

\bibitem{rao}
M. M. Rao, Conditional measure and applications, Marcel Dekker,
New York, 1993.


\bibitem{z}
A. C. Zaanen, Integration, 2nd ed., North-Holland, Amsterdam,
1967.
\end{thebibliography}
\end{document}